# GEOMETRIC MEAN METHOD FOR JUDGEMENT MATRICES: FORMULAS FOR ERRORS

I.L. Tomashevskii

*Institute of Mathematics, Information and Space Technologies,*

*Northern (Arctic) Federal University, Arkhangelsk 163000, Russia*

*E-mail:* tomashevil@gmail.com

## ABSTRACT

The geometric mean method (GMM) and the eigenvector method (EM) are well-known approaches to deriving information from pairwise comparison matrices in decision making processes. However, the original algorithms of these methods are logically incomplete and have significant drawbacks: their actual numerical errors are unknown and their reliability is doubted by different rank reversal phenomena including the GMM-EM ones.

Recently (Tomashevskii 2015: Eur J Oper Res 240:774-780), the actual numerical errors were found for EM, and it was shown that all EM rank reversal phenomena have the same cause and are eliminated when the errors are taken into account. In this paper the similar approach is used for GMM: we associate GMM with some standard measuring procedure, analyze this procedure from the viewpoint of measurement theory, and find the actual GMM errors. We show that the GMM and the GMM-EM rank reversal phenomena are eliminated when the GMM and the EM errors are taken into account.

The GMM decision support tool, which has all components of a standard measuring tool, is composed of pairwise comparisons as an initial measuring procedure, GMM as a data processor, and the obtained formulas for GMM errors as an error indicator. This tool is analogous to the EM decision support tool received earlier in the above paper. It is shown that the EM and the GMM versions of the tool are equally suitable to measure and rank any comparable elements with positive numerical values.

We also analyze the Geometric Consistency Index usually used in the original GMM to measure of the inconsistency and to accept or reject an inconsistent pairwise comparison matrix, and show that this index is not an acceptable GMM error indicator.

## 1. Introduction

Pairwise comparison matrices are handy to represent the preferences of experts in decision making processes. Different methods are applied to derive finished information from these matrices (e.g., see Golany & Kress 1993; Choo & Wedley 2004; Wang, Luo & Xu 2011; Kou & Lin 2013). They transform expert's judgments on the considered elements (alternatives) to the numerical values (priorities). In this paper, we are focused on two methods. One of them is the eigenvector method (EM) (Saaty 1980), which derives values (priorities) $\omega_1,..,\omega_n$ of comparable elements $\Omega_1,\ldots,\Omega_n$ as the solution of the eigenvalue problem

$$\sum_{k=1}^{n} a_{ik}\omega_k = \lambda_{\max}\omega_i\,,\quad i=1,\ldots,n, \tag{1}$$

for the corresponding pairwise comparison matrix $A=\| a_{ik} \|$ with the principle eigenvalue $\lambda_{\max}$. Another is the geometric mean method (GMM) also called the approximate eigenvector method (Saaty 1980; Crawford & Williams 1985; Kumar & Ganesh 1996) or the logarithmic least squares method. This method is handy for group decision making (e.g., see Ramanatham & Ganesh 1994; Forman & Peniwati 1998; Dong et al. 2010) and derives priorities of elements $\Omega_1,..,\Omega_n$ as

$$\omega_i = C\left(\prod_{k=1}^{n} a_{ik}\right)^{1/n},\quad i=1,\ldots,n\,, \tag{2}$$

where $C$ is a normalization factor.

Initially, the derived values $\omega_1,..,\omega_n$ are intended to estimate and rank the elements:

$$\omega_k > \omega_i \;\Rightarrow\; \Omega_k \succ \Omega_i\,,\quad \omega_k < \omega_i \;\Rightarrow\; \Omega_k \prec \Omega_i\,.$$

However, the problem is that the values $\omega_1,..,\omega_n$ are exact only when the pairwise comparison matrix $A=\|a_{ik}\|$ is perfectly consistent, i.e., the transitivity rule $a_{ik}a_{kr}=a_{ir}$ holds for all comparisons (Saaty 1977). Any inconsistency entails the errors $\Delta\omega_1,..,\Delta\omega_n$ of the values $\omega_1,..,\omega_n$. For this reason, the inequality $\omega_k > \omega_i$ (or $\omega_k < \omega_i$) implies the reliable ranking $\Omega_k \succ \Omega_i$ (or $\Omega_k \prec \Omega_i$) of the elements $\Omega_k$ and $\Omega_i$ only if the errors $\Delta\omega_k$ and $\Delta\omega_i$ are suffi-

ciently small, and the intervals $[\omega_k - \Delta\omega_k, \omega_k + \Delta\omega_k]$ and $[\omega_i - \Delta\omega_i, \omega_i + \Delta\omega_i]$ are disjoint, i.e., $\Delta\omega_k + \Delta\omega_i < |\omega_k - \omega_i|$. If this condition is strongly violated, for instance,

$$\omega_k \pm \Delta\omega_k$$

−−−−−−−−−−−(\\\\\\\\\\•\\\\\\\\\\)−−→

$$\omega_i \pm \Delta\omega_i$$

−−−−(\\\\\\\\\\\\\\•\\\\\\\\\\\\\\)−−−−→

then the inequality $\omega_k > \omega_i$ (or $\omega_k < \omega_i$) does not contain any reliable information on the ranking of the elements $\Omega_k$ and $\Omega_i$, and the ranking outside the probability theory is not possible.

As shown (Tomashevskii 2015), the invalid use of only the values $\omega_1,..,\omega_n$ for a ranking is the same cause of all EM rank reversal phenomena: the rank reversal problem for scale inversion or the right-left eigenvector asymmetry (Johnson et al. 1979; Saaty 1980), the rank reversal phenomenon caused by the addition or deletion of an element under consideration (Hochbaum et al. 2006; Raharjo et al. 2005), and the reversal of “order of intensity of preference” (Bana e Costa and Vansnick 2008). As shown, these phenomena are eliminated when the actual EM errors $\Delta\omega_1,..,\Delta\omega_n$:

$$\Delta\omega_i = \sqrt{\frac{1}{n-1}\sum_{k=1}^{n}\left(\frac{n}{\lambda_{\max}}a_{ik}\omega_k - \omega_i\right)^2}\ , \quad i = 1,..,n\,, \tag{3}$$

and the actual values $\omega_1 \pm \Delta\omega_1,..,\omega_n \pm \Delta\omega_n$ are taken into account. The formulas (3) complete the original EM algorithm (1) into the full decision support tool (Tomashevskii 2015), which has all properties of a standard measuring tool: it generates quantitative estimations $\omega_1,..,\omega_n$ and indicates their errors $\Delta\omega_1,..,\Delta\omega_n$. At present GMM has significant drawbacks, which are similar to the ones of the original EM.

Firstly, the actual numerical errors of the derived GMM values $\omega_1,..,\omega_n$ (2) is unknown. GMM uses the Geometric Consistency Index (*GCI*) (Crawford & Williams 1985; Aguarón & Moreno-Jiménez 2003) to accept or reject an inconsistent pairwise comparison matrix. The *GCI* is a direct analog of the Consistency Ratio (*CR*) (Saaty 1980) and is used similarly the Saaty’s criterion of $CR \le 0.1$. However, as proved in (Tomashevskii 2015), *CR* is not an ac-ceptable EM error indicator.[1] As a consequence, the *GCI* cannot be an acceptable GMM error indicator. In addition, the *GCI* is only some heuristic criterion, which does not intend to detect actual numerical errors of the derived values $\omega_1,..,\omega_n$.

Secondly, GMM as well as the original EM allows different rank reversal phenomena. Moreover, for the same pairwise comparison matrix, the ranking of the GMM values $\omega_1,..,\omega_n$ is not necessarily the same as the ranking of the EM ones. This phenomenon was demonstrated by Saaty and Vargas (1984), and Saaty (1990) and was the subject of some discussion (see e.g. Barzilai 1997; Lootsma 1999; Saaty 2005; Dijkstra 2013). In particular, Barzilai (1997) claims that "the geometric mean is the only method for deriving weights from multiplicative pairwise comparisons which satisfies fundamental consistency requirements", and Saaty (2005, sec.2-4) expresses their opinion that EM is the "only plausible candidate for representing priorities derived from a near consistent pairwise comparison matrix".

In this paper, we use a correspondence between GMM and some "matrix" measuring procedure, analyze this procedure from the viewpoint of measurement theory, and find actual numerical errors $\Delta\omega_1,..,\Delta\omega_n$ of the GMM values $\omega_1,..,\omega_n$. We show that GMM and EM are equally suitable to measure and rank any comparable elements with positive numerical values when their actual numerical errors $\Delta\omega_1,..,\Delta\omega_n$ are taken into account, and show that GMM can be considered as a full decision support tool.

## 2. Matrix measurements, GMM matrix measuring tool

Let $\Omega_1,\Omega_2,...,\Omega_n$ be comparable elements with positive numerical values on some common scale. Suppose we have the measuring procedure

$$\begin{Bmatrix} \Omega_1 \\ ... \\ \Omega_n \end{Bmatrix} \rightarrow \left\langle \begin{matrix} \text{measuring} \\ \text{procedure} \end{matrix} \right\rangle \rightarrow \begin{Bmatrix} \omega_1^{(1)},.., \omega_1^{(n)} \\ ... \\ \omega_n^{(1)},.., \omega_n^{(n)} \end{Bmatrix},$$

which realizes some matrix measurement from $n$ independent measurements of each elements. Then, according to the measurement theory (e.g., see Taylor 1997), we can use the mean

[1] The eligibility of the *CR* criterion has been much debated before (Monsuur 1997; Karapetrovic & Rosenbloom 1999; Kwiesielewicz & van Uden 2004; Bana e Costa & Vansnick 2008; Bozoki & Rapcsak 2008).

$$\overline{\omega}_i = \frac{1}{n}\sum_{k=1}^{n}\omega_i^{(k)} \tag{4}$$

as an approximate value of the element $\Omega_i$, the standard deviation

$$\Delta\overline{\omega}_i = \sqrt{\frac{1}{n-1}\sum_{k=1}^{n}\left(\omega_i^{(k)} - \overline{\omega}_i\right)^2} \tag{5}$$

as its error, and

$$\overline{\omega}_i \pm \Delta\overline{\omega}_i \quad , \quad i = 1,..,n, \tag{6}$$

as actual values of the elements $\Omega_1, \Omega_2, \ldots, \Omega_n$. If we modernize the measuring procedure:

$$\begin{Bmatrix} \Omega_1 \\ \ldots \\ \Omega_n \end{Bmatrix} \to \left\langle \begin{matrix} \text{ln - measuring} \\ \text{procedure} \end{matrix} \right\rangle \to \begin{Bmatrix} \ln\omega_1^{(1)},.., \ln\omega_1^{(n)} \\ \ldots \\ \ln\omega_n^{(1)},.., \ln\omega_n^{(n)} \end{Bmatrix},$$

then we can use the mean

$$\ln\omega_i^* = \frac{1}{n}\sum_{k=1}^{n}\ln\omega_i^{(k)} = \ln\left(\prod_{k=1}^{n}\omega_i^{(k)}\right)^{1/n}$$

with

$$\omega_i^* = \left(\prod_{k=1}^{n}\omega_i^{(k)}\right)^{1/n} \tag{7}$$

as an approximate logarithmic value of the element $\Omega_i$, the standard deviation

$$\Delta\ln\omega_i^* = \sqrt{\frac{1}{n-1}\sum_{k=1}^{n}\left(\ln\omega_i^{(k)} - \ln\omega_i^*\right)^2} = \sqrt{\frac{1}{n-1}\sum_{k=1}^{n}\ln^2\left(\omega_i^{(k)}/\omega_i^*\right)} \tag{8}$$

as its error, and

$$\ln\omega_i^* \pm \Delta\ln\omega_i^*, \quad i = 1,..,n,$$

as actual logarithmic values of the elements $\Omega_1,..,\Omega_n$. For actual normal values of the elements $\Omega_1,..,\Omega_n$ we obtain

$$\omega_i \pm \Delta\omega_i = \exp(\ln\omega_i^* \pm \Delta\ln\omega_i^*) = \omega_i^* \cdot \exp(\pm\Delta\ln\omega_i^*) \quad , \quad i = 1,..,n, \tag{9}$$

where

$$\omega_i = \omega_i^* \frac{\exp(\Delta\ln\omega_i^*) + \exp(-\Delta\ln\omega_i^*)}{2} \equiv \omega_i^* \cdot \mathrm{ch}(\Delta\ln\omega_i^*), \tag{10}$$

$$\Delta\omega_i = \omega_i^* \frac{\exp(\Delta\ln\omega_i^*) - \exp(-\Delta\ln\omega_i^*)}{2} \equiv \omega_i^* \cdot \mathrm{sh}(\Delta\ln\omega_i^*) \ . \tag{11}$$

(the actual values (6) and (9) are in close agreement for all $0 - 0.7$ relative errors).

The measuring procedure together with (10), (11) compose some measuring tool

$$\begin{Bmatrix} \Omega_1 \\ \dots \\ \Omega_n \end{Bmatrix} \to \left\langle \begin{matrix} \text{measuring} \\ \text{procedure} \end{matrix} \right\rangle \to \begin{Bmatrix} \omega_1^{(1)},.., \omega_1^{(n)} \\ \dots \\ \omega_n^{(1)},..,\omega_n^{(n)} \end{Bmatrix} \to \left\langle \begin{matrix} \text{data processor (10)} \\ \text{error indicator (11)} \end{matrix} \right\rangle \to \begin{Bmatrix} \omega_1 \pm \Delta\omega_1 \\ \dots \\ \omega_n \pm \Delta\omega_n \end{Bmatrix}. \quad (12)$$

Recently in (Tomashevskii 2015), the similar measuring tool with data processor in form (4) and the error indicator in form (5) was adapted by using EM to the measuring procedure in the form of pairwise comparisons. GMM can be used to adapt the measuring tool (12) to one. This possibility ensues from the following.

**Proposition.** *For any positive values* $\{\omega_i^{(k)}\}_{i,k=1}^n$ *with the geometric means* $\omega_1^*,..,\omega_n^*$ *defined by* (7) *and any C>0 there exists a unique* $n \times n$ *positive matrix* $A = \|a_{ik}\|$ *such that*

$$\omega_i^* = C\left(\prod_{k=1}^n a_{ik}\right)^{1/n}, \quad \omega_i^{(k)} = \frac{1}{\lambda} a_{ik}\omega_k^*, \quad i,k = 1,..,n, \quad \lambda = \left(\prod_{i,r=1}^n a_{ir}\right)^{1/n^2}. \quad (13)$$

*Conversely, for any* $n \times n$ *positive matrix* $A = \|a_{ik}\|$ *and any positive constant* $C$ *there exist unique positive values* $\{\omega_i^{(k)}\}_{i,k=1}^n$ *with the geometric means* $\omega_1^*,..,\omega_n^*$ *such that* (13) *holds* (see the proof in Appendix A).

This proposition reveals the correspondence between all measuring tools (12) and positive matrices $A = \| a_{ik} \|$ such that

1) by (7), (8), (13), the means $\omega_1,..,\omega_n$ (10) of $\Omega_1,..,\Omega_n$ and these errors (11) are equal to

$$\omega_i = \omega_i^* \operatorname{ch}\Delta_i \ , \quad (14)$$

$$\Delta\omega_i = \omega_i^* \operatorname{sh}\Delta_i \ , \quad (15)$$

$$\omega_i^* = C\left(\prod_{k=1}^n a_{ik}\right)^{1/n}, \quad \Delta_i = \sqrt{\frac{1}{n-1}\sum_{k=1}^n \ln^2\left(\frac{1}{\lambda} a_{ik}\frac{\omega_k^*}{\omega_i^*}\right)}, \quad (16)$$

for all $i = 1,...,n$ (for reciprocally symmetric matrices, $a_{ki} = 1/a_{ik}$, we have $\lambda \equiv 1$; the constant $C$ is normalization factor),

2) if a measuring procedure in (12) is absolutely precise (i.e., $\omega_1^{(k)} = \omega_1, \dots, \omega_n^{(k)} = \omega_n$ for all $k$), then (see (7),(13)) $\omega_1^* = \omega_1, \dots, \omega_n^* = \omega_n$ and

$$a_{ik} = \text{const} \cdot \frac{\omega_i}{\omega_k}, \quad i,k = 1,..,n, \tag{17}$$

i.e., $A$ is similar to a pairwise comparison matrix.

The above means that any practicable comparison procedure (algorithm) for comparable elements $\Omega_1,..,\Omega_n$, which generates a positive matrix $A = \|a_{ik}\|$ with functionally independent entries and makes possible the precise limit (17), along with (14), (15) composes the following "GMM matrix measuring tool"

$$\begin{Bmatrix} \Omega_1 \\ \dots \\ \Omega_n \end{Bmatrix} \to \left\langle \begin{matrix} \text{comparison} \\ \text{procedure} \end{matrix} \right\rangle \to \begin{Vmatrix} a_{11} \dots & a_{1n} \\ \dots & \\ a_{n1} \dots & a_{nn} \end{Vmatrix} \to \left\langle \begin{matrix} \text{data processor}\,(14) \\ \text{error indicator}\,(15) \end{matrix} \right\rangle \to \begin{Bmatrix} \omega_1 \pm \Delta\omega_1 \\ \dots \\ \omega_n \pm \Delta\omega_n \end{Bmatrix}. \tag{18}$$

This tool has all properties of a standard measuring tool: it generates quantitative estimations $\omega_1,..,\omega_n$ and indicates their errors $\Delta\omega_1,..,\Delta\omega_n$, and is an analog of the "EM matrix measuring tool" with a data processor in form (1) and an error indicator in form (3) received earlier in (Tomashevskii 2015).

## 3. The geometric mean method as a GMM matrix measuring tool

It is obvious that the standard pairwise comparison procedure

$$a_{ik} = \frac{\text{approximate value of } \Omega_i}{\text{approximate value of } \Omega_k} \tag{19}$$

used in decision making processes automatically satisfies to the precise limit (17) : if the procedure (19) is absolutely precise then the approximate value of $\Omega_i$ tends to the exact value $(= \omega_i)$ and

$$a_{ik} = \frac{\text{exact value of } \Omega_i}{\text{exact value of } \Omega_k} = \frac{\omega_i}{\omega_k}, \quad i,k = 1,..,n.$$

Moreover, for small inconsistency (small errors), we get (see (14)-(16)) $\Delta_i \approx 0$, $\text{ch}\Delta_i \approx 1$, $\text{sh}\Delta_i \approx \Delta_i$, and $\omega_i \approx \omega_i^*$, $\Delta\omega_i \approx \omega_i^* \Delta_i$. It means that the original GMM (2) is a data processor of the "small inconsistency" versions of the GMM matrix measuring tool (18). For a large inconsistency, the original GMM values in form (2) transforms into (14) and the formulas (15) indicate their errors.

### 3.1 On rank reversal problems

Our testing shows that the GMM and the EM matrix measuring tools generate the actual values $\omega_1 \pm \Delta\omega_1,..,\ \omega_n \pm \Delta\omega_n$ , which are in close agreement with each other for all permissible errors. This fact is interesting due to the GMM-EM rank reversal phenomenon (Saaty & Vargas 1984; Saaty 1990; Hovanov et al. 2008): using the same pairwise comparison matrix, the original EM and the original GMM provide different rankings (based only on the ranking of the mean values $\omega_1,..,\omega_n$). For illustration, we test the example (Saaty & Vargas 1984) with the reciprocal pairwise comparison matrix

$$A = \begin{pmatrix} 1 & 1/6 & 1/3 & 1/8 & 5 \\ \cdot & 1 & 2 & 1 & 8 \\ \cdot & \cdot & 1 & 1/2 & 5 \\ \cdot & \cdot & \cdot & 1 & 5 \\ \cdot & \cdot & \cdot & \cdot & 1 \end{pmatrix} \tag{20}$$

The normalized priorities ($\sum \omega_k = 1$) generated by the original GMM/EM and the GMM/EM matrix measuring tools are following:

$$\begin{array}{ccc} \text{original} & \text{GMM}(2) & \text{EM}(1) \\ \omega_1 & 0.073 & 0.081 \\ \omega_2 & 0.358 & 0.346 \\ \omega_3 & 0.187 & 0.180 \\ \omega_4 & 0.345 & 0.355 \\ \omega_5 & 0.036 & 0.038 \\ \multicolumn{3}{c}{\omega_2 \leftrightarrow \omega_4 \ \textit{rank reversal}} \end{array} \xrightarrow{\text{actual}} \begin{array}{cc} \text{GMM}(14)+(15) & \text{EM}(1)+(3) \\ 0.080 \pm 0.039 & 0.081 \pm 0.056 \\ 0.344 \pm 0.052 & 0.346 \pm 0.063 \\ 0.179 \pm 0.017 & 0.180 \pm 0.027 \\ 0.357 \pm 0.142 & 0.355 \pm 0.155 \\ 0.040 \pm 0.020 & 0.038 \pm 0.018 \\ \multicolumn{2}{c}{\textit{The close agreement}} \end{array} \tag{21}$$

We see that all the GMM and the EM actual results are in close agreement. The GMM-EM rank reversal problem is eliminated.

The original GMM as well as the original EM also allows rank reversal phenomena when the ranking of the initial elements is changed by the addition (or deletion) of some element (Hochbaum et al. 2006) and the reversal of "order of intensity of preference" (Bana e Costa & Vansnick 2008). It is shown (Tomashevskii 2015), in the case of EM, these phenomena have the same cause: the invalid use of only the principle eigenvector entries $\omega_1,..,\omega_n$ for ranking, and are self-eliminated when the actual values $\omega_1 \pm \Delta\omega_1,..,\omega_n \pm \Delta\omega_n$ are taken into account. Our testing shows that, in the case of GMM, the cause and the consequences are absolutely the same: all rank reversal phenomena arise when the errors $\Delta\omega_1,..,\Delta\omega_n$ are so large that the ranking using only the values $\omega_1,..,\omega_n$ is not valid.

### 3.2 On the Geometric Consistency Index as an error indicator

The GMM errors $\Delta\omega_1,..,\Delta\omega_n$ are a consequence of pairwise comparison matrix inconsistencies. For reliable ranking of the elements $\Omega_1,..,\Omega_n$, these errors (and the inconsistency) must be sufficiently small (see Introduction). Thus far, the original GMM (2) use the Geometric Consistency Index

$$GCI_n = \frac{2}{(n-1)(n-2)} \sum_{i=1}^{n} \sum_{k=i}^{n} \ln^2\left(a_{ik}\omega_k / \omega_i\right)$$

(Crawford & Williams 1985) to measure of the inconsistency and to accept or reject an inconsistent $n \times n$ pairwise comparison matrix. This heuristic index is a direct analog to the Saaty's Consistency Index (Saaty 1980) $CI_n = (\lambda_{\max} - n)/(n-1)$ used for EM. Saaty (1980) claimed that if $CI_n / RI_n \le 0.1$, where $RI_n$ is the average random Consistency Index derived from a sample of size 500 of randomly generated $n \times n$ reciprocal matrices, then the pairwise comparison matrix is acceptable. Using this criterion, Aguarón and Moreno-Jiménez (2003) determines a corresponding threshold for $GCI_n$: $GCI_3 \le 0.0314$, $GCI_4 \le 0.0352$, $GCI_n \le 0.037$ $(n > 4)$. The mathematical analyze of the Saaty's criterion shows that this criterion is not an acceptable EM error indicator (Tomashevskii 2015). As a consequence, the Aguarón and Moreno-Jiménez's criterion cannot be an acceptable GMM error indicator. For purposes of illustration, we return to the above example with the matrix (20). In the case under consideration, $CI_5 / RI_5 = 0.08 \le 0.1$ and $GCI_5 = 0.27 \le 0.037$ (according to Saaty, and Aguarón and Moreno-Jiménez, the matrix *A* is acceptable). However, according to (21), the errors $\Delta\omega_1,..,\Delta\omega_n$ of $\omega_1,..,\omega_5$ are such that the reliable ranking of the elements under consideration is not possible. It is clear from the graphical representation of the actual values (21):

```
      ω5±Δω5       ω3±Δω3              ω2±Δω2
 -|-(\\•\\)-----(\\•\\)------(\\\\\\•\\\\\\)-----------|-→
         ω1±Δω1                          ω4±Δω4
 -|--(\\\\•\\\\)------(\\\\\\\\\\\\\\\\\•\\\\\\\\\\\\\\\\\)--|-→
  0                                                   0.55
```

It means that the Aguarón and Moreno-Jiménez's criterion enables to accept the inconsistent pairwise comparison matrices *A*, which contains the unacceptable errors.

## 4. Another version of the GMM matrix measuring tool

GMM with a pairwise comparison matrix $A=\|a_{ik}\|$ (19) is one of versions of the matrix measuring tool (18). Another version is GMM with the transposed matrix $A^T=\|\tilde{a}_{ik}\|$,

$$\tilde{a}_{ik}=a_{ki}=\frac{\text{approximate value of } \Omega_k}{\text{approximate value of } \Omega_i}=\frac{1/\text{approximate value of } \Omega_i}{1/\text{approximate value of } \Omega_k}, \qquad (22)$$

which satisfies the precise limit (10) for the inverse values $\tilde{\omega}_1,..,\tilde{\omega}_n$ ($\tilde{\omega}_i \approx 1/$approximate value of $\Omega_i$) of the elements $\Omega_1,...,\Omega_n$ (see the comment to (19)).

According to (14)-(16), for the inverse values $\tilde{\omega}_1,..,\tilde{\omega}_n$ and their errors $\Delta\tilde{\omega}_1,..,\Delta\tilde{\omega}_n$ we get

$$\tilde{\omega}_i=\tilde{\omega}_i^*\,\text{ch}\tilde{\Delta}_i,\quad \Delta\tilde{\omega}_i=\tilde{\omega}_i^*\,\text{sh}\tilde{\Delta}_i\ ,\quad \tilde{\omega}_i^*=C\left(\prod_{k=1}^{n}\tilde{a}_{ik}\right)^{1/n},\quad \tilde{\Delta}_i=\sqrt{\frac{1}{n-1}\sum_{k=1}^{n}\ln^2\left(\frac{1}{\lambda}\tilde{a}_{ik}\frac{\tilde{\omega}_k^*}{\tilde{\omega}_i^*}\right)}\ .$$

In the case a reciprocal pairwise comparison matrix *A*, $a_{ki}=1/a_{ik}$, and *C*=1, we obtain (see (14)-(16), (22))

$$\tilde{a}_{ik}=a_{ki}=\frac{1}{a_{ik}},\quad \tilde{\omega}_i^*=\left(\prod_{k=1}^{n}\tilde{a}_{ik}\right)^{1/n}=\left(\prod_{k=1}^{n}a_{ik}\right)^{-1/n}=\frac{1}{\omega_i^*},\quad \tilde{\Delta}_i=\Delta_i,\quad i=1,...,n,$$

and

$$\tilde{\omega}_i=\frac{1}{\omega_i^*}\text{ch}\Delta_i,\quad \Delta\tilde{\omega}_i=\frac{1}{\omega_i^*}\text{sh}\Delta_i\,. \qquad (23)$$

For absolutely precise measurements, the matrix measuring tools with matrices (22) and (19) generate mutually inverse values of the compared elements $\Omega_1,..,\Omega_n$: if one tool generates the values $\tilde{\omega}_1,..,\tilde{\omega}_n$ then another tool generates the values $\omega_1,..,\omega_n$, such that $\omega_i=1/\tilde{\omega}_i$ for all *i*. In the case of an inconsistent pairwise comparison matrix *A*, the inverse values $\tilde{\omega}_i \pm \Delta\tilde{\omega}_i$ are transformed to the normal values

$$w_i \pm \Delta w_i=\frac{1}{\tilde{\omega}_i \mp \Delta\tilde{\omega}_i}=\frac{1}{\tilde{\omega}_i \mp \Delta\tilde{\omega}_i}\cdot\frac{\tilde{\omega}_i \pm \Delta\tilde{\omega}_i}{\tilde{\omega}_i \pm \Delta\tilde{\omega}_i}=\frac{\tilde{\omega}_i}{\tilde{\omega}_i^2-\Delta\tilde{\omega}_i^2}\pm\frac{\Delta\tilde{\omega}_i}{\tilde{\omega}_i^2-\Delta\tilde{\omega}_i^2}$$

For reciprocal pairwise comparison matrices, according to (23), (14), (15), we get

$$w_i \pm \Delta w_i=\frac{\omega_i^*\,\text{ch}\Delta_i}{\text{ch}^2\Delta_i-\text{sh}^2\Delta_i}\pm\frac{\omega_i^*\,\text{sh}\Delta_i}{\text{ch}^2\Delta_i-\text{sh}^2\Delta_i}=\omega_i^*\,\text{ch}\Delta_i \pm \omega_i^*\,\text{sh}\Delta_i=\omega_i \pm \Delta\omega_i\,.$$

Both versions of the GMM matrix measuring tool are equally suitable to measure and rank any comparable elements $\Omega_1,..,\Omega_n$ with positive numerical values.

## 5. Conclusion

Our goal is achieved. We showed that GMM as well as EM is a full tool to measure and rank any comparable elements with positive numerical values. In the form (18), this tool has all properties of a standard measuring tool and generates the values which are in close agreement with corresponding EM values. This tool can be generalized to group decision making.

## APPENDIX A

Proof Proposition. Let $C$, $\{\omega_i^{(k)}\}_{i,k=1}^n$ be any positive values and

$$\omega_i^* = \left(\prod_{k=1}^n \omega_i^{(k)}\right)^{1/n} .$$

Consider the positive values

$$\lambda = \frac{1}{C}\left(\prod_{k=1}^n \omega_k^*\right)^{1/n} , \quad a_{ik} = \lambda\omega_i^{(k)} / \omega_k^*, \quad i,k = 1,..,n. \qquad \text{(A.1)}$$

Clearly,

$$\left(\prod_{k=1}^n a_{ik}\right)^{1/n} = \lambda\omega_i^* / \left(\prod_{k=1}^n \omega_k^*\right)^{1/n} , \quad i = 1,..,n ,$$

and

$$\omega_i^* = C\left(\prod_{k=1}^n a_{ik}\right)^{1/n} . \qquad \text{(A.2)}$$

Moreover, from (A.1), (A.2) it is follows that

$$\lambda = \left(\prod_{i,k=1}^n a_{ik}\right)^{1/n^2} .$$

Thus, (13) holds.

Converselly, let $A = \|a_{ik}\|$ be a $n \times n$ positive matrix, $C$ be a positive constant and

$$\lambda = \left(\prod_{k,r=1}^n a_{kr}\right)^{1/n^2} , \quad \omega_i^{(k)} = \frac{C}{\lambda} a_{ik}\left(\prod_{r=1}^n a_{kr}\right)^{1/n} , \quad i,k = 1,..,n.$$

Then the geometric mean of $\left\{\omega_i^{(k)}\right\}_{k=1}^{n}$ is equal to

$$\omega_i^* = \left(\prod_{k=1}^{n} \omega_i^{(k)}\right)^{1/n} = \frac{C}{\lambda}\left(\prod_{k=1}^{n} a_{ik}\left(\prod_{r=1}^{n} a_{kr}\right)^{1/n}\right)^{1/n} = C\left(\prod_{k=1}^{n} a_{ik}\right)^{1/n}$$

Thus, (13) holds. This completes the proof of the proposition.